%% file: klein_bottles.tex
\documentclass[12pt,reqno]{article}
\include{head}

\title{Lagrangian Klein bottles in \(S^2\times S^2\)}
\author{Nikolas Adaloglou and Jonathan David Evans}

\newcommand{\pg}{\paragraph{\hspace{-0.6cm}}}
\newcommand{\Addresses}{{% additional braces for segregating \footnotesize
  \bigskip
  \footnotesize
  N. Adaloglou, \textsc{Mathematical Institute, Leiden University
  }, \texttt{n.adaloglou@math.leidenuniv.nl}\par\nopagebreak
  J.~D.~Evans, \textsc{University of Lancaster}, \texttt{j.d.evans@lancaster.ac.uk}

}}

\begin{document}

\maketitle

\begin{abstract}
  \noindent We use Luttinger surgery to show that there are no
  Lagrangian Klein bottles in \(S^2\times S^2\) in the
  \(\mathbb{Z}_2\)-homology class of an \(S^2\)-factor if the
  symplectic area of that factor is at least twice that of the
  other.
\end{abstract}

\section{Introduction}

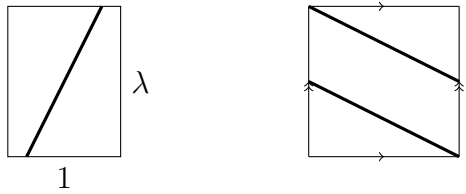
\begin{figure}[htb]
  \begin{center}
    \begin{tikzpicture}
      \draw (0,0) -- (1.5,0) node[midway,below] {\(1\)} --
      (1.5,2) node [midway,right] {\(\lambda\)} -- (0,2) -- cycle;
      \draw[very thick] (0.25,0) -- (1.25,2);
      \begin{scope}[shift={(4,0)}]
        \draw[->-] (0,0) -- (2,0);
        \draw[->-] (0,2) -- (2,2);
        \draw[->>-] (0,0) -- (0,2);
        \draw[->>-] (2,0) -- (2,2);
        \draw[very thick] (0,2) -- (2,1);
        \draw[very thick] (0,1) -- (2,0);
      \end{scope}
    \end{tikzpicture}
  \end{center}
  \caption{The base (left) and fibre (right) of the moment map
    \(\mu\colon X\to \RR^2\). A visible Lagrangian Klein bottle
    lives over the bold line in the base (slope \(2\)) and
    intersects each fibre in the bold line shown on the right
    (slope \(-1/2\)).}
  \label{fig:visible}
\end{figure}

\pg Let \(X=S^2\times S^2\) and let \(\omega_\lambda\) be the
product symplectic form which gives the factors areas \(1\) and
\(\lambda\) respectively. Define the homology classes
\begin{align*}
  \alpha &\coloneqq [S^2\times\{p\}]\\
  \beta &\coloneqq [\{p\}\times S^2].
\end{align*}
If \(\lambda < 2\) then there is a Lagrangian Klein bottle in the
homology class \(\beta\); one can construct this as a visible
Lagrangian submanifold \cite{EvansKB} (see Figure \ref{fig:visible}
below\footnote{The picture in {\cite[Fig. 2]{EvansKB}} is wrong and
should be rotated by \(90\) degrees otherwise the Klein bottle lives
in the \(\ZZ_2\)-homology class \(\alpha\).}) or in several other
ways \cite{CastroUrbano,DaiHoLi,Goodman}.

The line over which the visible Lagrangian Klein bottle
projects must have slope \(2\) and connect the bottom and top
edges, so it can be drawn if and only if \(\lambda <2\). For
this reason, the second author conjectured in \cite{EvansKB}
that there is no Lagrangian Klein bottle in the class \(\beta\)
when \(\lambda \geq 2\). We will prove this.

\paragraph{Theorem.}\label{thm:main} {\em If \(L\subset X\) is a
  Lagrangian Klein bottle for \(\omega_\lambda\) in the homology
  class \(\beta\) then \(\lambda<2\).}

\pg Indeed, there are no null-homologous Klein bottles in \(X\) by
Shevchishin \cite{Shevchishin} and Nemirovski \cite{Nemirovski}, and
there are none in the \(\ZZ_2\)-homology class \(\alpha+\beta\)
because they would violate the Audin identity \cite{Audin}:
\[\chi(L) = [L]\cdot[L]\mod 4\] which holds for
any totally real embedded submanifold \(L\) of an almost complex
surface. Here \(\chi\) is the Euler characteristic and
\([L]\cdot[L]\mod 4\) denotes the Pontryagin square of the
\(\ZZ_2\)-homology class. Therefore Theorem \ref{thm:main} gives a
complete picture of which homology classes are inhabited by Lagrangian
Klein bottles for which symplectic forms.

\pg There is also a complete understanding of which homology
classes are inhabited by {\em totally real} Klein bottles,
worked out by Derdzinski and Januszkiewicz {\cite[Proposition
  29.1]{DerdJan}}. If we allow immersions then all homology
classes can be represented, but there are restrictions on the
Maslov classes. If we allow only embeddings then the class
\(\alpha+\beta\) cannot be represented but the others can. The
restriction on Maslov classes will be vitally important to us,
so we review this in \ref{pg:maslov_constraint} once we have
established more notation.
    
\pg The usual toric diagrams for other Hirzebruch surfaces also
carry visible Klein bottles. We can perform a sequence of almost
toric mutations to get from such a diagram to either a rectangle
(for an even Hirzebruch surface) or a triangle with its corner
truncated (for an odd Hirzebruch surface), and this sequence of
almost toric mutations happens away from the visible Klein
bottles, so the families of visible Klein bottles coming from
different Hirzebruch surfaces are symplectomorphic to each
other. Note that in the case of odd Hirzebruch surfaces, the
visible Klein bottles are also \textit{real}, i.e. fixed loci of
anti-symplectic involutions.

\begin{figure}[htb]
  \begin{center}
    \begin{tikzpicture}
      \draw[->-] (0,0) -- (2,0) node [midway,below] {\(A\)};
      \draw[->-] (0,2) -- (2,2);
      \draw[->>-] (0,0) -- (0,2) node [midway,left] {\(B\)};
      \draw[->>-] (2,2) -- (2,0);
    \end{tikzpicture}
  \end{center}
  \caption{The Klein bottle as an identification space.}
  \label{fig:K}
\end{figure}
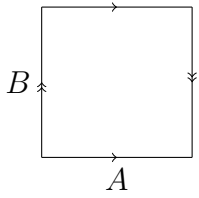

\paragraph{Klein bottles.} We will think of the Klein bottle as
the quotient of \(\RR^2\) (coordinates \((\varphi,\psi)\)) by
the action generated by the transformations
\[(\varphi,\psi)\mapsto (\varphi+1,-\psi),\qquad
  (\varphi,\psi)\mapsto (\varphi,\psi+1).\] We can represent
this quotient as a square with its sides identified as in Figure
\ref{fig:K}; we write \(A\) for the homology class of the loop
\(t\mapsto (t,0)\) and \(B\) for the homology class of the loop
\(t\mapsto (0,t)\). These generate the first homology and
satisfy \(2B=0\), \(A\cdot B=1\), \(B^2=0\), \(A^2=1\).

Because \(H_2(L;\ZZ)=H_1(X;\ZZ)=0\), the long exact sequence of
the pair \((X,L)\) splits off a short exact sequence
\[0\to H_2(X;\ZZ)=\ZZ^2\to H_2(X,L;\ZZ)\to
  H_1(L;\ZZ)=\ZZ\oplus\ZZ_2\to 0.\] In particular, for any class
\(H_1(L;\ZZ)\) there exists a class \(H_2(X,L;\ZZ)\) which maps
to it under the connecting homomorphism \(\partial\). The set of
such \(H_2(X,L;\ZZ)\)-classes form a torsor over \(H_2(X;\ZZ)\).

\paragraph{Maslov class constraint.}\label{pg:maslov_constraint}
If \(\Sigma\subset X\) is an oriented surface with boundary on a totally
real submanifold \(L\subset X\) then it has a well-defined Maslov
index (see \ref{maslov-par} for a review). In fact, the reduction
modulo 4 of this Maslov index depends only on the boundary \(\partial
\Sigma\in H_1(L;\ZZ)\). Derdzinski and Januszkiewicz \cite{DerdJan} call
this mod \(4\) index \(\bm{i}(\partial D)\).

\begin{proof}[Proof that \(\bm{i}(\partial \Sigma)\) is well-defined modulo \(4\)]
  Recall that \(\partial^{-1}([\partial
    \Sigma])\) is a torsor over \(H_2(X;\ZZ)\). All classes in
  \(H_2(X;\ZZ)\) have Maslov index equal to zero modulo \(4\), since the
  Maslov index on a class \(C\in H_2(X;\ZZ)\) agrees with
  \(2c_1(X)\cdot[C]\) and \(c_1(X)=2(\alpha+\beta)\).
\end{proof}

\pg We will be interested in surfaces with boundary homologous
to the loop \(B\). Derdzinski and Januszkiewicz show that for
totally real immersions of the Klein bottle into
\(S^2\times S^2\) in the \(\ZZ_2\)-classes \(0\) and
\(\alpha+\beta\), we must have \(\bm{i}(B)=0\mod 4\), but that
for totally real immersions in the \(\ZZ_2\)-classes \(\alpha\)
and \(\beta\), we have \(\bm{i}(B)=2\mod 4\). This comes from
the characterisation of the set \(\mathcal{Z}(\Sigma,M)\) of
possible Maslov/homology class pairs that appears after
{\cite[Theorem 2.2]{DerdJan}}, which is spelled out for Klein
bottles in \(S^2\times S^2\) in the paragraph before
{\cite[Proposition 29.1]{DerdJan}}. We reproduce the specific
parts of their argument we need in
\ref{pg:derdjan_0}--\ref{pg:derdjan_3} below for convenience.

\begin{figure}[htb]
  \begin{center}
    \begin{tikzpicture}
      \draw (0,0) -- (1.5,0) -- (1.5,2) -- (0,2) -- cycle; \draw[very thick] (0.25,0) -- (1.25,2) node [midway,left] {\(\ell\)};
      \begin{scope}[shift={(4,0)}]
        \draw[->-] (0,0) -- (2,0);
        \draw[->-] (0,2) -- (2,2);
        \draw[->>-] (0,0) -- (0,2);
        \draw[->>-] (2,0) -- (2,2);
        \draw[very thick] (0,2) -- (2,1);
        \draw[very thick] (0,1) -- (2,0);
        \draw[very thick,dashed] (1,0) -- (1,0.5);
        \draw[very thick,dotted] (1,0.5) -- (1,1.5);
        \draw[very thick,dashed] (1,1.5) -- (1,2);
      \end{scope}
    \end{tikzpicture}
  \end{center}
  \caption{A visible symplectic sphere living over \(\ell\) made
  up of two discs (dotted and dashed) with boundary on the
  visible Lagrangian \(L\).}
  \label{fig:visible_sphere}
\end{figure}
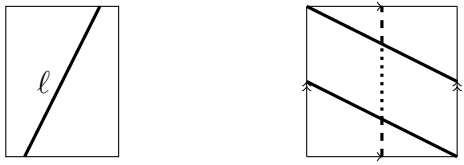

\paragraph{Example.}\label{exm:visible_bottle}
To see that \(\bm{i}(B)=2\mod 4\) for the visible Klein bottles,
consider the visible symplectic (but not holomorphic)
\(2\)-sphere living over the same line \(\ell\) as our visible
bottle which intersects every (regular) fibre in a vertical loop
as shown in Figure \ref{fig:visible_sphere}. This sphere lives
in the homology class \(\beta\) and intersects \(L\) in a circle
in the homology class \(B\), which separates the sphere into two
discs (one dotted, one dashed in Figure
\ref{fig:visible_sphere}). These discs form part of a
1-parameter family, as the vertical line moves horizontally; in
particular they have the same Maslov index. Since the sphere has
Chern number \(2\) (and hence Maslov index \(4\)), each disc has
Maslov index \(2\). Again, we emphasise these are symplectic but
not holomorphic discs.

\paragraph{Outline of the proof.}
We outline the proof of Theorem \ref{thm:main} here, leaving two
facts to be established later (Section
\ref{sct:remainder}). Recall that given a Lagrangian Klein
bottle \(L\subset X\) there is a symplectic surgery called
Luttinger surgery which excises a Weinstein neighbourhood of
\(L\) and reglues it with a twist. This surgery is based on the
surgery for Lagrangian tori introduced by Luttinger
\cite{Luttinger}, modified for Klein bottles by Nemirovski
\cite{Nemirovski}. We will write
\((\tilde{X},\tilde{\omega})\) for the result of this
Luttinger surgery.

We will first prove (\ref{prp:lutt_1}) that \(\tilde{X}\) is
diffeomorphic to the first Hirzebruch surface
\(\FF_1\cong\CC\PP^2\#\overline{\CC\PP}^2\). The second homology
of \(\FF_1\) is generated by two classes \(H\) and \(E\) with
squares \(1\) and \(-1\) respectively. We next prove
(\ref{prp:lutt_2}) that
\[\int_E\tilde{\omega}=1-\frac{\lambda}{2}.\] By a result of
Li and Liu {\cite[Theorem A]{LiLiuRuled}}, building on the
seminal work of Taubes (specifically {\cite[Proposition
  4.5]{Taubes}}), the fact that \(E\) is represented by a
smoothly embedded sphere and has \(E^2=-1\) and \(c_1(E) = 1\)
implies that the class \(E\) is represented by a symplectic
sphere. This then implies that \(\int_E\tilde{\omega}>0\),
that is \(1-\lambda/2>0\), or \(\lambda<2\), as required.

\section{Klein bottle Luttinger surgery}

In this section, we will review the Klein bottle Luttinger
surgery of Nemirovski \cite{Nemirovski} in a form closer to the
exposition of Auroux, Donaldson and Katzarkov \cite{ADK}. We
will establish some basic but important properties about how
Chern and Maslov classes are related before and after the
surgery.

\pg\label{pg:F_def_1} Let \(D_r\) denote the square
\([-r,r]\times[-r,r]\), with coordinates \((p_1,p_2)\) and let
\(T^2\) be the 2-torus with \(1\)-periodic coordinates
\((q_1,q_2)\). Consider the symplectic manifold
\(\left(D_{2r}\setminus D_r\right)\times T^2\) with symplectic
form \(\sum dp_i\wedge dq_i\) (this is a punctured subdomain in
\(T^*T^2\)). Let \(\chi\colon [-r,r]\to [0,1]\) be a smooth,
monotonically increasing function satisfying
\[\chi(p)=\begin{cases}
    0&\mbox{ if }p<-r/3\\
    1&\mbox{ if }p>r/3.
  \end{cases}\] To ensure compatibility with the Klein bottle
case later, we will further assume that \(\chi(p)-1/2\) is an
odd function, that is
\[\chi(-p)=-\chi(p)+1.\] As observed by Auroux, Donaldson and
Katzarkov \cite{ADK}, the diffeomorphism
\[F(p_1,p_2,q_1,q_2)=
  \begin{cases}
    (p_1,p_2,q_1,q_2+\chi(p_2))&\mbox{ if }p_1>0\\
  (p_1,p_2,q_1,q_2)&\mbox{ if }p_1<0\\
  \end{cases}
\] preserves the symplectic form:
\[F^*\left(\sum dp_i\wedge dq_i\right) = dp_1\wedge
  dq_1+dp_2\wedge (dq_2+\chi'(p_2)dp_2)=\sum dp_i\wedge dq_i.\]

\pg\label{pg:F_def_2} We will think of \(T^2\) as a double cover
of the Klein bottle \(L\), with deck transformation
\((q_1,q_2)\mapsto (q_1+1/2,-q_2)\). This deck transformation
induces a symplectomorphism of \(D_{2r}\times T^2\):
\[\Psi(p_1,p_2,q_1,q_2) = (p_1,-p_2,q_1+1/2,-q_2)\] which
commutes with \(F\). To see this, observe that
\begin{align*}
  F(\Psi(p_1,p_2,q_1,q_2))&=(p_1,-p_2,q_1+1/2,-q_2+\chi(-p_2))\\
  \Psi(F(p_1,p_2,q_1,q_2))&=(p_1,-p_2,q_1+1/2,-q_2-\chi(p_2)).
\end{align*}
Since by assumption \(\chi(-p_2)=-\chi(p_2)+1\), we have
\[-q_2+\chi(-p_2)=-q_2-\chi(p_2)+1,\] so the final coordinates
of \(F\circ\Psi\) and \(\Psi\circ F\) coincide (recall that they
are periodic modulo \(1\)). Let us write \(U_r\) for the
neighbourhood of the zero-section in \(T^*L\) given by
quotienting \(D_{r}\times T^2\) by \(\Psi\) and we continue to
write \(F\) for the symplectomorphism of \(U_r\) induced by
\(F\colon D_r\times T^2\to D_r\times T^2\).

\pg If \(X\) is a symplectic manifold containing a Lagrangian
Klein bottle \(L\) then we can find a symplectic embedding
(Weinstein neighbourhood) \(i\colon U_{2r}\to X\) for some
\(r\). Let \(\mathcal{U}=i(U_{2r})\), let
\(\mathcal{V}=X\setminus i(U_r)\), and let
\(\mathcal{W}=\mathcal{U}\cap \mathcal{V}\cong U_{2r}\setminus
U_r\). We can form the surgered manifold
\(\tilde{X}=\mathcal{U}\cup_{\mathcal{W}} \mathcal{V}\), where a
point \(x\in \mathcal{W}\) is identified with
\(x\in \mathcal{V}\) and with \(F(x)\in \mathcal{U}\).

\pg Consider the hypersurface
\(N = i(\{p_1^2+p_2^2=3r^2\}) \subset W\). Whilst the boundaries
of \(\mathcal{U}\) and \(\mathcal{V}\) are only piecewise
smooth, this is a smooth contact hypersurface isomorphic to the
radius \(r\) cosphere bundle of \(L\). The hypersurface \(N\)
separates \(X\) into two closed sets \(U\) (containing \(L\))
and \(V\) which are deformation retracts of \(\mathcal{U}\) and
\(\mathcal{V}\). We can think of \(V\) as a subset of
both \(X\) and \(\tilde{X}\) which is filled by two different
symplectic fillings of \(N\) which we call \(U\) (before
surgery) and \(\tilde{U}\) (after surgery).

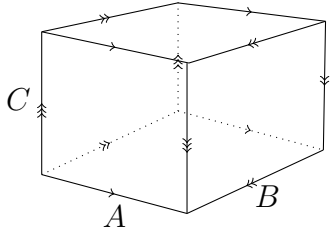
\begin{figure}[htb]
  \begin{center}
    \begin{tikzpicture}
      \draw[->-] (0,0) -- ++ (-15:2) node (a) {};
      \draw[dotted,->>-] (0,0) -- ++ (25:2) node (b) {};
      \draw[->>>-] (0,0) -- ++ (0,1.9) node (c) {};
      \draw[dotted,->-] (b.center) -- ++ (-15:2);
      \draw[dotted,->>>-] (b.center) -- ++ (0,1.5);
      \draw[-<<-] (a.center) -- ++ (25:2) node (g) {};
      \draw [-<<<-] (a.center) -- ++ (0,2);
      \draw[->-] (c.center) -- ++ (-12:2) node (d) {};
      \draw[->>-] (c.center) -- ++ (12:1.85) node (e) {};
      \draw[-<<-] (d.center) -- ++ (17:1.9) node (f) {};
      \draw[->>-] (f.center) -- (g.center);
      \draw[->-] (e.center) -- (f.center);
      \node at (-15:1) [below] {\(A\)};
      \node at (0,1) [left] {\(C\)};
      \node at (3,0) [below] {\(B\)};
    \end{tikzpicture}
  \end{center}
  % >>
  \caption{The space \(N\) as an identification space.}
  \label{fig:NL}
\end{figure}

\pg Note that \(N\) is a circle bundle over \(L\); we will use
coordinates \((q_1,q_2,\theta)\) on \(N\), where the coordinates
are understood modulo the identifications
\begin{align*}
  (q_1,q_2,\theta)&\mapsto (q_1+1/2,-q_2,-\theta)\\
  (q_1,q_2,\theta)&\mapsto (q_1,q_2+1,\theta)\\
  (q_1,q_2,\theta)&\mapsto (q_1,q_2,\theta+1)
\end{align*}
shown in Figure
\ref{fig:NL}.

We write \(A\), \(B\) and \(C\) for the homology classes of the three
loops \(t\mapsto (t/2,0,0)\), \(t\mapsto (0,t,0)\) and \(t\mapsto
(0,0,t)\).

\paragraph{Lemma.} {\em From the cell decomposition inherited
  from the cube (respectively the square), we can easily
  compute:
  
  \begin{center}
    \begin{tabular}{c||c|c}
      & \(N\) & \(L\) \\
      \hline\hline
      \(H_3(-;\ZZ)\) & \(\ZZ\) & \(0\)\\
      \(H_2(-;\ZZ)\) & \(\ZZ\) & \(0\)\\
      \(H_1(-;\ZZ)\) &  \(\langle A,B,C\,\mid\,2B=2C=0\rangle\) & \(\langle A,B\,\mid\,2B=0\rangle\)\\
      \hline
      \(H^3(-;\ZZ)\) & \(\ZZ\) & \(0\)\\
      \(H^2(-;\ZZ)\) & \(\ZZ\oplus\ZZ_2^2\) & \(\ZZ_2\)\\
      \(H^1(-;\ZZ)\) & \(\ZZ\) & \(\ZZ\)
    \end{tabular}
  \end{center}

  and the inclusion \(i\colon N\to U\simeq L\) induces the map
  \(A\mapsto A\), \(B\mapsto B\), \(C\mapsto 0\) on
  homology. Over \(\RR\), the pullback \(i^*\) induces an
  isomorphism on \(H^1\) and the zero map on \(H^2\). }

\paragraph{Lemma.} {\em Now suppose that
  \(H^1_{dR}(X)=H^3_{dR}(X)=0\). We have \(H^1_{dR}(V)\cong 0\)
  and \(H^2_{dR}(V)\cong H^2_{dR}(X)\oplus\RR\). Moreover, the
  pullback map \(H^2_{dR}(X)\to H^2_{dR}(V)\) is injective.}
\begin{proof}
  For the purposes of the proof it is more convenient to work
  with \(\mathcal{U}\) and \(\mathcal{V}\) than \(U\) and \(V\),
  but since they are deformation retracts the conclusions are
  the same. Since \(N\) is a deformation retract of
  \(\mathcal{U}\cap \mathcal{V}\), the relevant part of the
  Mayer--Vietoris sequence for \(X=\mathcal{U}\cup \mathcal{V}\)
  is:
  \begin{align}
    \nonumber\cdots\to H^1_{dR}(X)&\to H^1_{dR}(\mathcal{U})\oplus H^1_{dR}(\mathcal{V})\to H^1_{dR}(N)\to\\
    \nonumber\to H^2_{dR}(X)&\to H^2_{dR}(\mathcal{U})\oplus H^2_{dR}(\mathcal{V})\to H^2_{dR}(N)\to \\
    \label{eq:mv}\to H^3_{dR}(X)&\to\cdots
  \end{align}
  Using the facts that \(H^1_{dR}(X)=H_{dR}^3(X)=0\),
  \(H^2_{dR}(\mathcal{U})=0\), and that the map
  \(H^1_{dR}(\mathcal{U})\to H^1_{dR}(N)\) is an isomorphism, we deduce
  that \(H^1_{dR}(\mathcal{V})=0\) and split off a short exact sequence
  \[0\to H^2_{dR}(X)\to H^2_{dR}(\mathcal{V})\to H^2_{dR}(N)=\RR\to 0,\]
  which implies what is claimed.
\end{proof}

\paragraph{Corollary.}\label{pg:cohom} {\em The De Rham
  cohomology class of a 2-form on \(X\) is determined by its
  integrals over \(2\)-cycles in \(V\).}

\paragraph{Lemma.} \label{lma:conclusions} {\em Let
  \(L\subset X\) be a Lagrangian Klein bottle and let
  \(\tilde{X}\) be the result of performing Luttinger surgery
  along \(L\). Then
  \(H^1_{dR}(\tilde{X})=H^3_{dR}(\tilde{X})=0\),
  \(H^2_{dR}(\tilde{X})\cong H^2_{dR}(X)\), and the cohomology
  class of a 2-form on \(\tilde{X}\) is determined by its
  integrals over \(2\)-cycles in \(V\).}
\begin{proof}
  Write \(\tilde{X}\) as \(\tilde{U}\cup V\), where
  \(\tilde{U}\) is the Weinstein neighbourhood of the Klein
  bottle after surgery. The map \(N\to \tilde{U}\) still induces
  an isomorphism on \(H^1_{dR}\) and the maps
  \(H^k_{dR}(V)\to H^k_{dR}(N)\) are unchanged, so the
  Mayer--Vietoris sequence \eqref{eq:mv} (with \(\tilde{X}\)
  instead of \(X\) and \(\tilde{U}\) instead of \(U\)) gives the
  conclusions.
\end{proof}

\paragraph{Lemma.} \label{lma:nemirovski} {\em We have
  \(c_1(X)\cdot [\omega_\lambda] = c_1(\tilde{X})\cdot
  [\tilde{\omega}]\).} (This observation is due to Nemirovski
{\cite[Section 2.3]{Nemirovski}}; we have given details here for
convenience.)
\begin{proof}
  Represent \(c_1\) in each case by a Chern form
  \(\rho,\tilde{\rho}\) given by the curvature of an Hermitian
  line bundle. By picking the Hermitian metrics to agree over
  \(V\), we can ensure that \(\rho|_V=\tilde{\rho}|_V\). We also
  know that \(\omega_\lambda|_V=\tilde{\omega}|_V\), so that
  \[\int_C\omega_\lambda=\int_C\tilde{\omega}\qquad\mbox{and}\qquad
    \int_C\rho=\int_C\tilde{\rho}\] for any \(2\)-cycle \(C\subset
    V\). Since \(H^2_{dR}(U)=0\), both \(\omega_\lambda\) and \(\rho\)
    are cohomologous to 2-forms \(\omega'_\lambda\) and \(\rho'\)
    which vanish on \(X\setminus V\). Their restrictions to \(V\)
    extend (by zero) over \(\tilde{X}\setminus V\); let us write
    \(\tilde{\omega}'\) and \(\tilde{\rho}'\) for these
    extensions. Then, for any \(2\)-cycle \(C\subset V\), we have
  \[\int_C\tilde{\omega} = \int_C\omega_\lambda = \int_C\omega'_\lambda =
    \int_C\tilde{\omega}',\] and similarly
  \(\int_C\tilde{\rho}=\int_C\tilde{\rho}'\). By
  \ref{pg:cohom}, this shows that
  \([\tilde{\omega}]=[\tilde{\omega}']\) and
  \([\tilde{\rho}]=[\tilde{\rho}']\), so
  \[c_1(X)\cdot[\omega_\lambda]=\int_X\rho\wedge \omega_\lambda =
    \int_V\rho'\wedge \omega'_\lambda =
    \int_{\tilde{X}}\tilde{\rho}'\wedge\tilde{\omega}' =
    \int_{\tilde{X}}\tilde{\rho}\wedge\tilde{\omega} =
    c_1(\tilde{X})\cdot[\tilde{\omega}].\qedhere\]
\end{proof}

\paragraph{Corollary.}\label{cor:rational} {\em If
  \(c_1(X)\cdot [\omega]>0\) and \(H^1_{dR}(X)=0\) then
  \(\tilde{X}\) is a rational symplectic 4-manifold.}
\begin{proof}
  By Lemma \ref{lma:nemirovski},
  \(c_1(\tilde{X})\cdot[\tilde{\omega}]>0\), so by the
  Liu--Ohta--Ono Theorem ({\cite[Theorem B]{Liu}},
  {\cite[Theorem 1.2]{OhtaOno}}) the 4-manifold \(\tilde{X}\) is
  either rational or ruled. But by Lemma \ref{lma:conclusions},
  \(H^1_{dR}(\tilde{X})=0\), so \(\tilde{X}\) cannot be (a
  blow-up of) an irrational ruled surface. Therefore it is
  rational.
\end{proof}

\paragraph{Chern and Maslov classes.} Let \(X\) be an
almost complex manifold of real dimension \(2n\). Consider the
complex line bundle
\(V\coloneqq \left(\Lambda^n_\CC TX\right)^{\otimes 2}\) and let
\(E\subset V\) be its unit circle bundle. The first Chern class
of \(E\to X\) is \(2c_1(X)\), and if \(S\subset X\) is a
closed, oriented surface then \(2c_1(X)\cdot [S]\) is the
obstruction to finding a section of \(E|_S\); in other
words, \(2c_1(X)\cdot [S]\) is the signed count of zeros of
a generic section of \(V|_S\).

\pg\label{maslov-par}
If we have a surface \(\Sigma\) with boundary and a
nowhere-vanishing section \(\sigma\) of \(V\) defined along the
boundary \(\partial\Sigma\) then the Maslov index
\(\mu(\Sigma,\sigma)\) of the pair \((\Sigma,\sigma)\) is the
number of zeros of an extension of \(\sigma\) over
\(\Sigma\). One way of obtaining a nowhere-vanishing section of
\(V\) over a subset \(Y\subset X\) is from a field of unoriented
totally real \(n\)-planes \(\zeta\) on \(Y\): namely, we pick a
section of \(\left(\Lambda^n_\RR \zeta\right)^{\otimes 2}\) by
squaring a local orientation of \(\zeta\). Since \(\zeta\) is
totally real, this section is nowhere-vanishing when considered
as a section of \(V\).

For example, if \(\Sigma\) has boundary on a totally real
submanifold \(L\) then there is a {\em
  canonical}\footnote{Canonical up to a positive real scalar.}
nowhere-vanishing section \(\sigma_{can}\) defined along
\(\partial \Sigma\) coming from the field of tangent planes
\(TL\). The Maslov index \(\mu(\Sigma)\) is then defined to be
the Maslov index of the pair \(\mu(\Sigma,\sigma_{can})\).

\pg\label{pg:maslov_2} If two boundary conditions \(\sigma_0\)
and \(\sigma_1\) are homotopic through sections of
\(E|_{\partial \Sigma}\) then
\(\mu(\Sigma,\sigma_0)=\mu(\Sigma,\sigma_1)\). If \(L\) is a
Lagrangian submanifold and we choose a Weinstein neighbourhood
then one obvious choice for an alternative boundary condition
equivalent to \(\sigma_{can}\) is the vertical distribution of
tangent spaces to cotangent fibres. This actually gives us a
section of \(E\) over the whole of the Weinstein
neighbourhood. We call this section \(\xi_{can}\).

\paragraph{Lemma.}\label{lma:F_action} {\em Let \(F\) be the
  diffeomorphism of a punctured Weinstein neighbourhood of a
  Lagrangian Klein bottle defined in
  \ref{pg:F_def_1}--\ref{pg:F_def_2}. The field \(\xi_{can}\) is
  homotopic through sections of \(E\) to \(F_*\xi_{can}\).}
\begin{proof}
  The vertical distribution on the cotangent bundle is spanned
  by \(\partial_{p_1}\) and \(\partial_{p_2}\). We have
  \(F_*\partial_{p_1}=\partial_{p_1}\) and
  \(F_*\partial_{p_2}=\partial_{p_2}+\chi'(p_2)\partial_{q_2}\). The
  homotopy between these is given by taking \(\partial_{p_1}\)
  and \(\partial_{p_2}+t\chi'(p_2)\partial_{q_2}\), which are
  complex linearly-independent for all \(t\).
\end{proof}

\paragraph{Lemma.}\label{lma:key} {\em Let \(X\) be a symplectic
  \(4\)-manifold and \(L\subset X\) a Lagrangian Klein
  bottle. Let \(U\) be a closed Weinstein neighbourhood of \(L\)
  with boundary \(N\) and let \(V=\overline{X\setminus U}\). Let
  \(\Sigma\subset X\) be a smooth, oriented surface with
  boundary on \(L\); by making a small perturbation, assume that
  \(\Sigma\) is transverse to \(N\). Let
  \(\Sigma'=V\cap \Sigma\); this is a smooth surface with
  boundary on \(N\). Suppose that \(F(\partial\Sigma)\) is
  nullhomologous in \(\tilde{U}\), pick an oriented smooth
  surface \(T\subset\tilde{U}\) with
  \(\partial T=-F(\partial\Sigma)\), and let
  \(\tilde{S}=\Sigma'\cup T\). Then}
  \[2c_1(\tilde{X})\cdot[\tilde{S}]=\mu(\Sigma).\]
\begin{proof}
  Using \(\xi_{can}\) on \(U\) and on \(\tilde{U}\) produces
  homotopic nowhere-vanishing boundary conditions for
  \(\Sigma'\), regardless of whether the rest of the surface is
  \(\Sigma\cap U\) or \(\tilde{S}\cap
  \tilde{U}\). Therefore the obstruction to extending these
  boundary conditions is the same; in the one case it is
  \(\mu(\Sigma)\) and in the other it is
  \(2c_1(\tilde{X})\cdot[\tilde{S}]\).
\end{proof}

\pg\label{pg:derdjan_0} For the convenience of the reader, we now
give a summary of how Derdzinski and Januszkiewicz prove that a
surface \(\Sigma\) in \(S^2\times S^2\) with boundary on a
Lagrangian Klein bottle \(L\) in the \(\ZZ_2\)-homology class
\(\alpha\) or \(\beta\) must have Maslov index
\(\mu(\Sigma)=2\mod 4\).

\paragraph{Lemma.}\label{pg:derdjan_1} {\em For \(X=S^2\times S^2\) we have
  \(\pi_1(E)=\ZZ/4\).}
\begin{proof}
  This follows from the homotopy long exact sequence of the
  circle bundle \(E\to X\):
  \[0\to \pi_2(E)\to \pi_2(X)\to \ZZ\to\pi_1(E)\to 0.\] Under
  the identification \(H^2(X;\ZZ)\cong \OP{Hom}(\pi_2(X),\ZZ)\),
  the map \(\pi_2(X)=\ZZ^2\to \pi_1(S^1)=\ZZ\) is the first
  Chern class of the circle bundle \(E\), that is
  \(2c_1(X)=4(\alpha+\beta)\). Therefore its image is the
  subgroup \(4\ZZ\subset\ZZ\) and its cokernel is
  \(\pi_1(E)=\ZZ/4\).\qedhere
\end{proof}

\paragraph{Lemma.}\label{pg:derdjan_2} {\em Let \(\Sigma\) be an
  oriented surface with boundary on a totally real surface
  \(L\subset X=S^2\times S^2\). Let \(\ell\) be the loop in
  \(E\) given by restricting the section \(\xi_{can}\) to
  \(\partial\Sigma\). Then
  \(\mu(\Sigma)\mod 4=[\ell]\in\pi_1(E)\).}
\begin{proof}
  The bundles \(E|_{\partial\Sigma}\) and \(E|_\Sigma\) are
  trivial, and only one of the trivialisations of
  \(E|_{\partial\Sigma}\) is compatible with a trivialisation of
  \(E|_\Sigma\). The Maslov index \(\mu(\Sigma)\) can be
  interpreted as the winding number of the boundary condition
  \(\xi_{can}\) with respect to this trivialisation, that is the
  projection of \(\xi_{can}|_{\partial\Sigma}\) to the second
  factor in
  \(\pi_1(E|_{\partial\Sigma}) = \pi_1(\partial\Sigma)\times
  \ZZ\). The inclusion map \(\partial \Sigma\to E\) induces a
  map \(\pi_1(\partial\Sigma)\times\ZZ\to \pi_1(E)=\ZZ/4\),
  which coincides with reduction modulo \(4\) on the second
  factor, and \(\xi_{can}|_{\partial\Sigma}\) maps to
  \([\ell]\).
\end{proof}

\paragraph{Lemma.}\label{pg:derdjan_3} {\em Let
  \(L\subset X=S^2\times S^2\) be a totally real Klein bottle
  (or, more generally, nonorientable surface with
  \(\chi(L)=0\mod 4\)) in the \(\ZZ_2\)-homology class
  \(\alpha\) or \(\beta\), let \(B\subset L\) be the meridian
  loop (or, more generally, a 2-sided simple closed curve
  representing the unique torsion class in \(H_1(L;\ZZ)\)) and
  let \(\Sigma\subset X\) be a surface with
  \(\partial \Sigma=B\). Then \(\mu(\Sigma)=2\mod 4\).}
\begin{proof}
  Let \(\ell\subset E\) be the loop given by the section
  \(\xi_{can}|_B\).  Since \(B\) is 2-sided, the section
  \(\xi_{can}|_B\) admits a square root (that is, a lift to
  \(\Lambda^n_\CC TX\)), which means that
  \([\ell]\in\{0,2\}\subset\pi_1(E)=\ZZ/4\). Suppose that
  \([\ell]=0\in\pi_1(E)\); we will derive a contradiction. We can cut
  open \(L\) along \(B\) to obtain an {\em orientable} surface \(L'\)
  with two circular boundary components; let \(L''\) be the abstract
  orientable surface obtained by capping these circles off with
  discs. Since \([\ell]=0\in\pi_1(E)\), the canonical section
  \(\sigma_L|_{L'}\colon L'\to E\) extends to a map \(L''\to E\) which
  sends the two capping discs to a single disc \(\Delta\) which is a
  nullhomotopy of \(\ell\). We can then project \(L''\) to \(X\) to
  obtain a continuous map \(f\colon L''\to X\). Since \(f^*E\) admits
  a section (by construction) this means that \(c_1(f^*E)=0\). Since
  \(c_1(f^*E)=2c_1(X)\cdot f_*[L'']\) and \(2c_1(X)=4(\alpha+\beta)\),
  this tells us that \(f_*[L'']=k(\alpha-\beta)\) for some
  \(k\in\ZZ\), which implies that \(f_*[L'']\mod 2\) is either \(0\)
  or \(\alpha+\beta\). But since the two capping discs project to the
  same disc (with opposite orientations), \(f_*[L'']=[L]\mod 2\),
  which contradicts the fact that \([L]\) is either \(\alpha\mod 2\)
  or \(\beta\mod 2\).
\end{proof}

\section{Finishing the proof}
\label{sct:remainder}
We now provide the details that were missing from the sketch proof.

\paragraph{Proposition.}\label{prp:lutt_1} {\em Let
  \(X=S^2\times S^2\) with the symplectic form
  \(\omega_{\lambda}\) and suppose that \(L\subset X\) is a
  Lagrangian Klein bottle in the \(\ZZ_2\)-homology class
  \(\beta\). Let \(\tilde{X}\) be the result of Luttinger
  surgery. Then \(\tilde{X}\) is diffeomorphic to
  \(\FF_1=\CC\PP^2\#\overline{\CC\PP}^2\).}
\begin{proof}
  By Corollary \ref{cor:rational}, \(\tilde{X}\) is a rational
  symplectic 4-manifold and by Lemma \ref{lma:conclusions}, it
  has \(b_2=2\). Therefore \(\tilde{X}\) is diffeomorphic to
  either \(S^2\times S^2\) or to
  \(\FF_1=\CC\PP^2\#\overline{\CC\PP}^2\). We will exhibit a
  surface \(\tilde{\Sigma}\) in \(\tilde{X}\) with odd Chern
  number, which is only possible for \(\FF_1\) since
  \(c_1(S^2\times S^2)\) is divisible by \(2\) in integral
  cohomology.

  To find \(\tilde{\Sigma}\), we start with a surface \(D\) with
  boundary on \(L\) such that \([\partial D]=B\) and
  \(\mu(D)=2\mod 4\). It is possible that \(D\) intersects \(L\)
  at interior points of \(D\); if so, we perturb \(D\) so that
  it intersects \(L\) transversely on its interior and open out
  these intersections to form small nullhomologous boundary
  components on \(L\). We continue to call the resulting surface
  \(D\). Since the boundary of \(D\) still lies in the class
  \(B\) (because the punctures gave only nullhomologous boundary
  components), the intersection \(D\cap N\) lies either in the
  class \(B+C\) or the class \(B\). In the former case, we
  define \(\Sigma=D\). In the latter case, we take
  \(\Sigma=S\cup D\) where \(S\) is a sphere in the homology
  class \(\alpha\), perturbed to be transverse to \(L\), again
  with these intersections opened out into boundary components
  on \(L\). Since \(S\) must intersect \(L\) an odd number of
  times, this has the effect of changing the intersection
  \(N\cap (D\cup S)\) so that it lives in the class \(B+C\).

  Now \(F_*(B+C)=C\), which is nullhomologous in \(\tilde{U}\),
  so we can apply Lemma \ref{lma:key}: note that
  \(\mu(\Sigma)=\mu(D)+k\mu(\alpha)\) for some \(k\), but
  \(\mu(\alpha)=4\), so \(\mu(\Sigma)=\mu(D)\mod 4\) and by
  Lemma \ref{pg:derdjan_3}, \(\mu(D)=2\mod 4\). This implies
  that \(2c_1(\tilde{\Sigma})=2\mod 4\), so
  \(c_1(\tilde{\Sigma})\) is odd.
\end{proof}

\paragraph{Proposition.}\label{prp:lutt_2} {\em Let \(H\) and
  \(E\) be the classes in \(\FF_1\) coming from a line in
  \(\CC\PP^2\) and the exceptional curve. If \(L\subset X\) is a
  Lagrangian Klein bottle in the \(\ZZ_2\)-class
  \(\beta\) in \(X=S^2\times S^2\) with symplectic form
  \(\omega_{\lambda}\) then the Luttinger surgery
  \((\tilde{X},\tilde{\omega})\) satisfies
  \[\int_E\tilde{\omega}=1-\frac{\lambda}{2}.\]}
\begin{proof}
  The rational homology of \(V\) fits into a Mayer--Vietoris
  sequence
  \[H_2(N;\ZZ)\to H_2(U;\ZZ)\oplus H_2(V;\ZZ)\to H_2(X;\ZZ)\to
    H_1(N;\ZZ)\to H_1(U;\ZZ).\] Since \([L]=\beta\mod 2\), a
  generic sphere in the class \(\beta\) intersects \(L\)
  transversely an even number of times. Thus
  \(\partial \beta\in H_1(N;\ZZ)\) is an even multiple of \(C\);
  since \(2C=0\), this means that \(\beta\) lives in the image
  of \(H_2(V;\ZZ)\). Similarly, \(2\alpha\) lies in the image of
  \(H_2(V;\ZZ)\). In particular, this means that
  \(2\alpha\) and \(\beta\) can be represented respectively by
  surfaces \(S_{2\alpha}\) and \(S_\beta\) in \(V\) with
  \begin{alignat*}{3}
    S_{2\alpha}^2=S_\beta^2=0, &\qquad& c_1(X)\cdot
    S_{2\alpha}=4 &\qquad &\int_{S_{2\alpha}}\omega_\lambda = 2,\\
    S_{2\alpha}\cdot S_\beta=2, &\qquad & c_1(X)\cdot S_\beta
    =2,&\qquad & \int_{S_\beta}\omega_\lambda = \lambda.
  \end{alignat*}
  After surgery, since these surfaces were in \(V\), they can be
  thought of as surfaces \(\tilde{S}_{2\alpha}\) and
  \(\tilde{S}_\beta\) in \(\tilde{X}\). We still have
  \begin{alignat*}{3}
    \tilde{S}_{2\alpha}^2=\tilde{S}_\beta^2=0, &\qquad&
    c_1(\tilde{X})\cdot \tilde{S}_{2\alpha}=4, &\qquad&
    \int_{\tilde{S}_{2\alpha}}\tilde{\omega} = 2,\\
    \tilde{S}_{2\alpha}\cdot \tilde{S}_\beta=2, &\qquad&
    c_1(\tilde{X})\cdot \tilde{S}_\beta =2, &\qquad& 
    \int_{\tilde{S}_\beta}\tilde{\omega} = \lambda.
  \end{alignat*}
  Suppose that \([S_\beta]=aH+bE\). Then
  \(S_\beta^2=a^2-b^2=0\), so \(a=\pm b\), and
  \(c_1(\tilde{X})\cdot S_\beta = 3a\pm a=2\). Since
  \(a\in\ZZ\), we must have \(a=-b=1\), so \([S_\beta]=H-E\).

  Suppose that \(S_{2\alpha}=cH+dE\). In the same way, from
  \(S_{2\alpha}^2=0\) we have
  \(c=\pm d\), and from \(c_1(\tilde{X})\cdot[S_{2\alpha}]=4\)
  we have \(3c\pm d=4\), so either \(c=d=1\) or \(c=-d=2\), that
  is \(S_{2\alpha}=H+E\) or \(2H-2E\). Since \(S_\alpha\cdot
  S_\beta=2\), we must have \(S_{2\alpha}=H+E\).

  But now
  \[\int_{S_{2\alpha}}\tilde{\omega}=2= \int_H\tilde{\omega}+ \int_E\tilde{\omega}
    \qquad\mbox{and}\qquad \int_{S_{\beta}}\tilde{\omega}=
    \lambda= \int_H\tilde{\omega}-\int_E\tilde{\omega},\] so
  \[2-\lambda=2\int_E\tilde{\omega}\qquad\mbox{or}\qquad
    \int_E\tilde{\omega}= 1-\frac{\lambda}{2}.\qedhere\]
\end{proof}

\section{Acknowledgements}

N.A. is grateful to his PhD advisor Federica Pasquotto, for her
overall guidance and support; to Jo\'e Brendel, Marco Golla,
George Politopoulos and Johannes Hauber for helpful and
insightful discussions; and finally to J.E. for inviting him to
Lancaster in May 2023, where this project was
conceived. J.E. would like to thank Stefan Nemirovski for vital
conversations about totally real submanifolds and Maslov
classes. J.E. is supported by EPSRC Grant EP/W015749/1.

\bibliographystyle{plain}
\bibliography{klein_bottles}
\Addresses
\end{document}

%% file: head.tex
\usepackage{amssymb}
\usepackage{enumerate}
\usepackage{bm}
\usepackage{lscape}
\usepackage{booktabs,caption}
\usepackage{multirow}
\usepackage{graphicx}
\usepackage{amsmath,amsthm}
\usepackage{a4wide}
\usepackage{version}
\usepackage{parskip}
\usepackage{euscript}
\usepackage{mathrsfs}
\usepackage{tikz}
\usepackage{mathtools}
\usetikzlibrary{arrows,decorations.markings, matrix}
\usepackage{hyperref}
\usepackage{comment}

\usepackage[all]{xy}

\numberwithin{equation}{section}

\newcommand{\CC}{\mathbb{C}}
\newcommand{\FF}{\mathbb{F}}

\newcommand{\PP}{\mathbb{P}}

\newcommand{\RR}{\mathbb{R}}

\newcommand{\ZZ}{\mathbb{Z}}

\newcommand{\OP}{\operatorname}

\tikzset{->-/.style={decoration={
              markings,
              mark=at position .5 with {\arrow{>}}},postaction={decorate}}}

\tikzset{-<-/.style={decoration={
              markings,
              mark=at position .5 with {\arrow{<}}},postaction={decorate}}}

\tikzset{->>-/.style={decoration={
               markings,
               mark=at position .5 with {\arrow{>>}}},postaction={decorate}}}

\tikzset{-<<-/.style={decoration={
               markings,
               mark=at position .5 with {\arrow{<<}}},postaction={decorate}}}

\tikzset{->>>-/.style={decoration={
               markings,
               mark=at position .5 with {\arrow{>>>}}},postaction={decorate}}}

\tikzset{-<<<-/.style={decoration={
               markings,
               mark=at position .5 with {\arrow{<<<}}},postaction={decorate}}}